% plain TeX; Gowers and Maurey, July 1994.

\catcode`\@=11

\font\tenmsx=msam10
\font\sevenmsx=msam7
\font\fivemsx=msam5
\font\tenmsy=msbm10
\font\sevenmsy=msbm7
\font\fivemsy=msbm5
\newfam\msxfam
\newfam\msyfam
\textfont\msxfam=\tenmsx  \scriptfont\msxfam=\sevenmsx
  \scriptscriptfont\msxfam=\fivemsx
\textfont\msyfam=\tenmsy  \scriptfont\msyfam=\sevenmsy
  \scriptscriptfont\msyfam=\fivemsy

\def\hexnumber@#1{\ifcase#1 0\or1\or2\or3\or4\or5\or6\or7\or8\or9\or
	A\or B\or C\or D\or E\or F\fi }

\edef\msx@{\hexnumber@\msxfam}
\edef\msy@{\hexnumber@\msyfam}

\mathchardef\square="0\msx@03
\mathchardef\leqslant="3\msx@36
\mathchardef\geqslant="3\msx@3E
\def\Bbb{\ifmmode\let\next\Bbb@\else
 \def\next{\errmessage{Use \string\Bbb\space only in math mode}}\fi\next}
\def\Bbb@#1{{\Bbb@@{#1}}}
\def\Bbb@@#1{\fam\msyfam#1}

\catcode`\@=12

\magnification=\magstep1
\hsize=6.5truein
\hoffset=0.0truein
\baselineskip 1.4\normalbaselineskip

\tolerance=10000
\def\sqr{$\vcenter{\hrule height .3mm
\hbox {\vrule width .3mm height 2mm \kern 2mm
\vrule width .3mm} \hrule height .3mm}$}
 
\def \Proof{\noindent {\bf Proof.\ \ }}
\def \R{{\Bbb R}}
\def \N{{\Bbb N}}
\def \C{{\Bbb C}}
\def \T{{\Bbb T}}
\def \Z{{\Bbb Z}}
\def \bfa{{\bf a}}
\def \bfe{{\bf e}}
\def \a{\alpha}
\def \b{\beta}
\def \g{\gamma}
\def \d{\delta}
\def \e{\epsilon}
\def\la{\lambda}

\def \s{\sigma}
\def \G{\Gamma}
\def\D{\Delta}
\def \ca{{\cal A}}
\def\ce{{\cal E}}
\def \ci{{\cal I}}
\def \cl {{\cal L}}
\def\co{{\cal O}}
\def \cs{{\cal S}}
\def\ct{{\cal T}}
\def \cu{{\cal U}}

\def \nm#1{\left\|#1\right\|}
\def \Bnm#1{\Bigl\|#1\Bigr\|}
\def \bgnm#1{\biggl\|#1\biggr\|}
\def \nmm#1{\left|\!\left|\!\left|#1\right|\!\right|\!\right|}
\def\nrm#1#2{\|#1\|_{(#2)}}
\def\bnrm#1#2{\left\|#1\right\|_{(#2)}}

\def \sq #1{(#1_n)_{n=1}^\infty}
\def \sqf #1 #2 #3{(#1_#2)_{#2=1}^#3}
\def \seq#1#2{#1_1,\dots,#1_#2}
\def \sleq#1#2{#1_1<\dots<#1_#2}
\def \speq#1#2{#1_1+\dots+#1_#2}
\def \sm#1#2{\sum_{#1=1}^#2}

\def \ge{\geqslant}
\def \le{\leqslant}
\def \ra{\rightarrow}

\def \ux{{\underline x}}

\def \codim{\mathop{\rm codim}}
\def\conv{\mathop{\rm conv}\nolimits}
\def \diag{\mathop{\rm diag}}

\def\flc{\mathop{\rightarrow}\limits}
\def\fl#1{\flc^{#1}}
\def \ind{\mathop{\rm ind}} 
\def \ker{\mathop{\rm ker}}
\def\mod{\mathop{\rm mod}}
\def\ran{\mathop{\rm ran}\nolimits}
\def\supp{\mathop{\rm supp}\nolimits}

\centerline {\bf BANACH SPACES WITH SMALL SPACES OF OPERATORS}
\medskip
\centerline {W. T. GOWERS and B. MAUREY}
\vskip 1truecm

\noindent {\bf Abstract.} {\sl For a certain class of algebras
$\cal A$ we give a method for constructing Banach spaces $X$ such that
every operator on $X$ is close to an operator in $\cal A$. This is
used to produce spaces with a small amount of structure. We present
several applications.
Amongst them are constructions of a new prime Banach
space, a space isomorphic to its subspaces of codimension two but not
to its hyperplanes and a space isomorphic to its cube but not to its
square.}
\vskip 2truecm

\noindent{\bf \S 1. Introduction}

This paper is a continuation of [GM], in which a space $X$ was
constructed with the property that every operator from a subspace $Y$
to $X$ was of the form $\la i+S$, where $i$ is the inclusion map and
$S$ is strictly singular. Among the easy consequences of this fact are
that $X$ contains no unconditional basic sequence, and, more
generally, that $X$ is hereditarily indecomposable: that is, no
subspace of $X$ admits any non-trivial projection.

One could say then that $X$ has as few operators as possible, which is
why it is a counterexample to many questions about general Banach
spaces.  However, there are other questions which assume some
structure for the space and then ask whether further structure
follows. For example, if $X$ has an unconditional basis, must it be
isomorphic to some proper subspace?  Various {\it ad hoc} techniques
have been developed by the authors for finding counterexamples to some
of these questions. The purpose of this paper is to be more
systematic.  We shall present a generalization of the main result of
[GM], which, roughly speaking, states that given an algebra of maps
satisfying certain conditions, one can replace the multiple of the
inclusion map in the statement above by the restriction to $Y$ of some
element of the algebra. This generalization has several
applications. Amongst them are constructions of a new prime Banach
space, a space isomorphic to its subspaces of codimension two but not
to its hyperplanes and a space isomorphic to its cube but not to its
square. A related argument shows that all the operators on the space
constructed in [G1], which has an unconditional basis, are of the form
$D+S$, where $D$ is diagonal and $S$ strictly singular. Note that the
diagonal operators must be continuous on such a space.

The paper will be organized as follows. In the next section we shall
introduce notation and some basic lemmas. These are similar to [GM]
but for technical reasons it was necessary to alter certain
definitions and prove statements that were not quite the same. One
result is that this paper is basically self-contained. In the third
section we state and prove the main result of the paper.  The
applications will be given in the fourth section. The fifth and final
section contains a discussion of the space mentioned above with an
unconditional basis. We also ask some questions about the possibility
of removing some of the restrictions needed for our main result.  To
understand the applications of our main theorem, it is not necessary
to understand the proof of the theorem, or even the definitions of the
spaces guaranteed to exist by it. The reader only interested in the
applications will be ready for them after reading the beginning of
Section 3 and perhaps skimming very lightly over Section 2.
\bigskip

\noindent {\bf \S 2. Notation and background}
\medskip

Let $c_{00}$ be the vector space of all (real or complex depending on
the context) sequences of finite support. Let $(\bfe_n)_{n=1}^\infty$
be the standard basis of $c_{00}$.  Given a vector $\bfa=\sm n \infty
a_n\bfe_n$ its {\sl support}, denoted $\supp(\bfa)$, is the set of $n$
such that $a_n\ne 0$. Given two subsets $E,F\subset\N$, we say that
$E<F$ if every element of $E$ is less than every element of $F$. If
$x,y\in c_{00}$, we say that $x<y$ if $\supp(x)<\supp(y)$. If $\sleq x
n$, then we say that the vectors $\seq x n$ are {\sl successive}. This
definition also applies to infinite sequences in an obvious way. An
infinite sequence of successive non-zero vectors is also called a {\sl
block basis} and a subspace generated by a block basis is a {\sl block
subspace}.  Given a subset $E\subset\N$ and a vector $\bfa$ as above,
we write $E\bfa$ for the vector $\sum_{n\in E}a_n\bfe_n$. That is, $E$
also stands for the coordinate projection associated with the set.  An
{\sl interval} of integers is a set of the form $\{n, n+1,\dots,m\}$
and the {\sl range} of a vector $x$, written $\ran(x)$, is the
smallest interval containing $\supp(x)$.

The following collection of functions was introduced by Schlumprecht
[S] (except for the technical condition (vi) below) and will be useful
here. It is the set $\cal F$ of functions $f:[1,\infty)\ra [1,\infty)$
satisfying the following conditions.

\itemitem {(i)} $f(1)=1$ and $f(x)<x$ for every $x>1$;
\itemitem {(ii)} $f$ is strictly increasing and tends to infinity;
\itemitem {(iii)} $\lim_{x\ra\infty}x^{-q}f(x)=0$ for every $q>0$;
\itemitem {(iv)} the function $x/f(x)$ is concave and non-decreasing;
\itemitem {(v)} $f(xy)\le f(x)f(y)$ for every $x,y\ge 1$;
\itemitem {(vi)} the right derivative of $f$ at $1$ is positive.
 
It is easy to check that $f(x)=\log_2(x+1)$ satisfies these
conditions, as does the function $\sqrt{f(x)}$. Note also that some of
the conditions above are redundant. In particular, it follows from the
other conditions that $f(x)$ and $x/f(x)$ are strictly increasing.
 
Now let $\cal X$ stand for the set of normed spaces $(c_{00},\nm .)$
such that the sequence $(\bfe_n)_{n=1}^\infty$ is a normalized
bimonotone basis. Given $X\in\cal X$ and $f\in\cal F$, we shall say
that $X$ {\sl satisfies a lower} $f$-{\sl estimate} if, given any
vector $x\in X$ and any sequence of intervals $\sleq E n$, we have
$\nm x\ge f(n)^{-1}\sm i n\nm{E_ix}$. Equivalently, if $\sleq x n$
then $\bgnm{\sm i n x_i}\ge f(n)^{-1}\sm i n\nm{x_i}$.

For $X\in \cal X$, $x\in X$ and every integer $N\ge 1$, let $$ \nrm xN
= \sup \sum_{i=1}^N \|E_i x\|,$$ where the supremum is extended to all
sequences $E_1,\ldots,E_N$ of successive intervals. Notice that, if
$X\in\cal X$, $x\in X$, $n\in\N$ and $E$ is an interval, then
$\nrm {Ex} n\le\nrm x n$.

For $0<\e\le 1$ and $f\in\cal F$, we say that a sequence
$x_1,\ldots,x_N$ of successive vectors {\sl satisfies the RIS($\e$)
condition for the function $f$} if there is a sequence
$(2N/f'(1))f^{-1}(N^2/\e^2) < n_1 <\cdots<n_N$ of integers (where
$f'(1)$ is the right derivative) such that $\nrm {x_i}{n_i} \le 1$ for
each $i=1,\ldots,N$ and $$\e {\sqrt {f(n_i)}} > |\ran(
\sum_{j=1}^{i-1} x_j)|$$ for every $i=2\ldots,N$. Observe that when
$x_1,\ldots,x_N$ satisfies the RIS($\e$) condition for some $f\in\cal
F$, then $Ex_1,\ldots,Ex_N$ also does for every interval $E$. When the
function $f$ is clear from the context, we shall simply say that $\seq
x N$ satisfies the RIS($\e$) condition.

Given $g\in\cal F$, $M\in\N$ and $X\in\cal X$, an $(M,g)$-{\sl form
on} $X$ is defined to be a functional $x^*$ of norm at most one which
can be written as $\sm j M x_j^*$ for a sequence $x_1^*<\dots<x_M^*$
of successive functionals all of which have norm at most $g(M)^{-1}$.
Observe that if $x^*$ is an $(M,g)$-form then $|x^*(x)|\le
g(M)^{-1}\nrm x M$ for any $x$. 

\proclaim Lemma 1. Let $f,g\in\cal F$ be such that $\sqrt f\le g$.
Assume that $x_1,\ldots,x_N$ satisfies the RIS($\e$) condition for
$f$. If $x^*$ is a $(k,g)$-form for some integer $k\ge 2$ then
$$ |x^*(\sum_{i=1}^N x_i)| \le \e + 1 + N/\sqrt{f(k)}.$$
In particular, $|x^*(x_1+\cdots+x_N)|< 1+2\e$ when
$k> f^{-1}(N^2/\e^2)$.

\Proof   Let $i$ be such that $n_i < k \le n_{i+1}$. Then,
since $\|x_j\|_{c_0} \le 1$ for every $j=1,\ldots,N$, we get
$$ | (x^*, \sum_{j=1}^{i-1} x_j )| \le \|x^*\|_\infty
\, |\supp(\sum_{j=1}^{i-1} x_j)| \le \e\sqrt{f(n_i)}/g(k) \le
\e,$$
$|(x^*, x_i)| \le 1$ and for every $j>i$
$$ |(x^*, x_j)| \le \nrm{x_j}{k} /g(k) \le 1/\sqrt{f(k)}.\eqno \square$$

\proclaim Lemma 2. Let $n\in\N$ and let $x\in X$ be a vector such
that $\nrm xn\le 1$. There exists a (non-negative) measure $w$ on
$A=\ran x$ such that $w(A)=1$ and such that $w(E)\ge\nm{Ex}$ for every
interval $E\subset A$ with $\nm{Ex}\ge n^{-1}$.

\Proof  Define $w^*(E)$ for any subinterval $E\subset A$ to be 
$\|Ex\|$ if $\|Ex\| \ge n^{-1}$, and $0$ otherwise.
It is enough to find a measure $w\ge w^*$ with $w(A)=1$.

We consider the linear programming problem of minimizing $w(A)$
subject to the family of constraints $ w(E) \ge w^*(E)$ for every
non-empty subinterval $E$ of $A$.  Let $\overline w$ be an optimal
solution for this problem and let
$${\cal J} = \{E; \ \overline w(E) = w^*(E)\}$$
be the set of active constraints for $\overline w$. It is a classical
fact that $A$ belongs to the closed convex cone generated by
the active constraints (identifying sets with characteristic
functions). We therefore have
$$A=\sum_{E\in{\cal J}}c_EE\ ,$$
with $c_E \ge 0$. Applying the difference operator $\D
f(x)=f(x)-f(x-1)$ to both sides of the above equation, we find that if
$E\in\cal J$, $c_E>0$ and $\max E<\max A$ then there exists $F\in\cal
J$ such that $c_F>0$ and $\min F=1+\max E$.  It follows by induction
that there exist $E_1<\ldots<E_l$ in $\cal J$ such that $A =
\bigcup_{i=1}^l E_i$.  

Since $\|x\|_{(n)} \le 1$, there are at most $n$ intervals
$E_i$ such that $w^*(E_i) > 0$, or equivalently such that
$\|E_i x\| \ge n^{-1}$. It follows that
$$ \overline w(A) = \sum_{i=1}^l \overline w(E_i) =
\sum_{i=1}^l w^*(E_i) \le \|x\|_{(n)} \le 1.$$
\hfill $\square$
\bigskip

\proclaim Lemma 3. Let $f,g\in\cal F$, $\sqrt f\le g$, and let
$\seq x N$ satisfy the RIS($\e$) condition for $f$. Let 
$x=\sm i Nx_i$ and suppose that
$$\nm{Ex}\le 1\vee\sup\bigl\{|x^*(Ex)|:
x^*\ \hbox{is a $(k,g)$-form}, k\ge 2\bigr\}$$
for every interval $E$. Then $\nm x\le (1+2\e)Ng(N)^{-1}$.

\Proof  Define $G(t)$ to be $t/g(t)$ when
$t\ge 1$ and $t$ when $t\le 1$. It is easy to check that $G$ is
concave, increasing and satisfies $G(st)\ge G(s)G(t)$ on the whole of
$\R_+$. 

Let $n_1$ be the first integer appearing in the RIS condition. We know
that $Nn_1^{-1} < 1$ by the RIS condition.  We have $\|x\|_{(n_1)} \le
N$ by the triangle inequality, thus we may find by Lemma 2 a measure
$w$ on $\ran(x)$ such that $w(\ran(x)) = N$ and $\|Ex \| \le w(E)$ for
every interval $E$ such that $\|Ex\| \ge Nn_1^{-1}$.  We call $w(E)$
the {\sl weight} of $E$.  We shall now show that $\nm{Ex}\le
(1+2\epsilon) G(w(E))$ for every interval $E$ such that $\|Ex\| \ge N
n_1^{-1}$.

If $ N n_1^{-1} \le \|Ex\| \le 1$, we have
$$ \|Ex\| = G(\|Ex\|) \le G(w(E))$$
and so the result holds in this case. Suppose that $E$ is a minimal
interval such that $\|Ex\| \ge N n_1^{-1}$ for which the inequality
fails.  We certainly have $\|Ex\| >1$, and by assumption $\|Ex\| >
(1+2\epsilon) G(w(E)) > (1+2\epsilon)$.  We therefore have a
$(k,g)$-form $x^*$ such that $(1+2\epsilon) G(w(E)) < |x^*(Ex)|$.  By
Lemma 1 and the definition of the RIS($\e$) condition, $k\le
f^{-1}(N^2/\e^2) \le f'(1) n_1/(2N)$.

It follows that we can find $\sleq E k$ with $\bigcup E_i=E$ and
$|x^*(Ex)| \le g(k)^{-1}\sm i k\nm{E_ix}$, by the definition of a
$(k,g)$-form. Let $w_i=w(E_i)$ for each $i$ and let $w=w(E)$. Since
$k\ge 2$ we may assume that no $E_i$ is equal to $E$. For each $i$ we
have either $\|E_ix\| \le N n_1^{-1}$ or, by the minimality of $E$,
$\nm{E_ix}\le (1+2\epsilon) G(w_i)$.  Let $A$ be the set of $i$ with
the first property and $B$ the complement of $A$. Let $s$ be the
cardinality of $A$.

By Jensen's inequality, we have
$$\sum_{i\in B}\nm{E_ix} \le (1+2\epsilon)
\sum_{i\in B}G(w_i) \le (1+2\epsilon) (k-s)G(w/(k-s))\ .
\eqno(*)$$
Therefore, setting $t=s/k$ and using the lower bound on $n_1$, we have
$$\eqalign{|x^*(Ex)| & \le (1+2\epsilon)
(1-s/k)G(k)G(w/(k-s))+sNn_1^{-1}\cr
& \le (1+2\epsilon) \big( (1-t)G(w/(1-t))+tf'(1)/2\big)\ .\cr}$$
(Observe that $s<k$, as otherwise
$|x^*(Ex)| \le kNn_1^{-1} < 1$).
Now note that $0<f'(1)/2\le g'(1)=G(1)-G'(1)$. Since $w\ge 1$ and
$0\le t<1$ it follows easily from the concavity of $G$ that
$|x^*(Ex)| \le (1+2\epsilon) G(w)$,
contradicting our assumption about the interval
$E$. The result follows. \hfill $\square$ \bigskip

\proclaim Lemma 4. Let $X\in\cal X$, satisfying a lower $f$-estimate
for some $f\in{\cal F}$. Then for every $n\in\N$ and $\e>0$,
every block subspace of $X$ contains a vector $x$ of finite support
such that $\nm x=1$ and $\nrm xn\le1+\e$. Hence, for every $N\in\N$,
every block subspace contains a sequence $\seq x N$ satisfying the
RIS($\e$) condition with $\nm{x_i}\ge (1+\e)^{-1}$.

\Proof  Without loss of generality $\e\le 1$. Let $m\ge 6n/\e$ be an 
integer. By a straightforward adaptation of
Lemma 3 of [GM] to a general $f\in{\cal F}$,
every block subspace contains a vector $x$ of norm 1 which
can be written as a sum $\speq x m$ of successive vectors where every
$x_i$ has norm at most $m^{-1}(1+\e/3)$. Let $\sleq E n$ be any
sequence of intervals whose union contains the support of $x$ and for
each $j\le n$ let $A_j=\{i:\supp(x_i)\subset E_j\}$ and let
$B_j=\{i:E_jx_i\ne 0\}$. By the triangle inequality and since the
basis of $X$ is bimonotone, $\nm{E_jx}\le\nm{\sum_{i\in
B_j}x_i}\le(1+\e/3)m^{-1}(|A_j|+2)$. Since $\sm j n|A_j|\le m$ we find
that $\sm j n\nm{E_jx}\le(1+\e/3)(1+2n/m)\le 1+\e$.
\hfill $\square$ \bigskip

\noindent {\bf \S 3. The main result}
\smallskip

We begin by defining a class of spaces, the adaptations which will
interest us of the space constructed in [GM].

Given two infinite sets $A,B\subset\N$, define the {\sl spread from}
$A$ {\sl to} $B$ to be the map on $c_{00}$ defined as follows. Let the
elements of $A$ and $B$ be written in increasing order respectively as
$\{a_1,a_2,\dots\}$ and $\{b_1,b_2,\dots\}$. Then $\bfe_n$ maps to
zero if $n\notin A$, and $\bfe_{a_k}$ maps to $\bfe_{b_k}$ for every
$k\in\N$. Denote this map by $S_{A,B}$. Let $P_A$ be the map
$S_{A,A}$, which is just the projection on to $A$. Note that
$S_{B,C}S_{A,B}=S_{A,C}$ and so $S_{B,A}S_{A,B}=P_A$. Note also that
$S_{B,A}$ is (formally) the adjoint of $S_{A,B}$.

Given any set $\cal S$ of spreads, we shall say that it is a {\sl
proper set} if it is closed under composition (note that this applies
to all compositions and not just those of the form $S_{B,C}S_{A,B}$)
and taking adjoints, and if, for every $(i,j)\ne(k,l)$, there are only
finitely many spreads $S\in\cal S$ for which $\bfe_i^*(S\bfe_j)\ne 0$ and
$\bfe_k^*(S\bfe_l)\ne 0$. A good example of such a set is the collection of
all spreads $S_{A,B}$ where $A=\{m,m+1,m+2,\dots\}$ and
$B=\{n,n+1,n+2,\dots\}$ for some $m,n\in\N$. This is the proper set
generated by the shift operator.

Given a Banach space $X$ satisfying a lower $f$-estimate for some
$f\in\cal F$, and given a subspace $Y\subset X$ generated by a block
basis, we will be interested in a seminorm $\nmm.$ defined on $L(Y,X)$
as follows. Let $\cl(Y)$ be the set of sequences $\sq x$ of successive
vectors in $Y$ such that $\nrm {x_n}{n}\le 1$. Now let $\nmm
T=\sup_{{\bf x}\in\cl(Y)}\lim\sup_n\nm{Tx_n}$.  The spaces we shall
consider satisfy lower $f$-estimates. Hence, by Lemma 4, all their
subspaces contain sequences in $\cl$ with norms bounded below (by 1/2
say). In such a space, if $\nmm T<\e$, then every subspace contains a
vector $x$ such that $\nm{Tx}<2\e\nm x$.  In particular, if $\nmm
T=0$, then $T$ is strictly singular.

The next theorem is the main one of the paper. For convenience, we now
fix $f\in\cal F$ for the rest of the paper to be the function
$f(x)=\log_2(x+1)$, as in the statement of the theorem.

\noindent {\bf Theorem 5.} {\sl Let $\cal S$ be a proper set of 
spreads. There exists a Banach space $X=X(\cal S)$ (satisfying a lower
$f$-estimate where $f(x)=\log_2(x+1)$) with the following three
properties.

(i) For every $x\in X$ and every $S_{A,B}\in\cal S$,
$\nm{S_{A,B}x}\le\nm x$, (and therefore $\nm{S_{A,B}x}=\nm x$ if
$\supp(x)\subset A$);

(ii) If $Y$ is a subspace of $X$ generated by a block basis, then
every operator from $Y$ to $X$ is in the $\nmm.$-closure of the set of
restrictions to $Y$ of operators in the algebra $\cal A$ generated by
$\cal S$. In particular, all operators on $X$ are
$\nmm.$-perturbations of operators in $\cal A$.

(iii) The seminorm $\nmm.$ satisfies the algebra inequality
$\nmm{UV}\le\nmm U\nmm V$.}
\bigskip

Notice a straightforward consequence of this result. If we write $\cal
G$ for the $\nmm.$-completion of $\cal A$ (after quotienting by
operators with $\nmm.$ zero) then $\cal G$ is a Banach algebra.  Given
$T\in L(X)$, we can find by (ii) a $\nmm.$-Cauchy sequence $\sq T$ of
operators in $\cal A$ such that $\nmm{T-T_n}\ra 0$. Let $\phi(T)$ be
the limit of $\sq T$ in $\cal G$. This map is clearly well-defined. It
follows easily from (iii) that it is also a unital algebra
homomorphism. The kernel of $\phi$ is the set of $T$ such that $\nmm
T=0$. The restriction of $\phi$ to $\ca$ is the identity (or more
accurately the embedding of $\ca$ into $\cal G$).  If $\cal A$ is
small, then, since the kernel of $\phi$ consists of small operators,
$L(X)$ is also small.

The first step in the proof of the theorem is to define the space
$X(\cal S)$.  First, we recall the definition from [GM] of the {\sl
special functionals} on a space $X\in\cal X$.  Let ${\bf Q}\subset
c_{00}$ be the set of sequences with rational coordinates and maximum
at most 1 in modulus.  Let $J\subset\N$ be a set such that, if $m<n$
and $m,n\in J$, then $\log\log\log n\ge 2m$.  Let us write $J$ in
increasing order as $\{j_1,j_2,\dots\}$. We shall also need
$f(j_1)>256$, where $f(x)$ is still the function $\log_2(x+1)$.  Now
let $K,L\subset J$ be the sets $\{j_1,j_3,j_5,\dots\}$ and
$\{j_2,j_4,j_6,\dots\}$.
 
Let $\sigma$ be an injection from the collection of finite sequences
of successive elements of ${\bf Q}$ to $L$.  Given $X\in\cal X$
such that $X$ satisfies a lower $f$-estimate and given an
integer $m\in\N$, let $A_m^*(X)$ be the set of functionals of the form
$f(m)^{-1}\sum_{i=1}^mx_i^*$ such that $x_1^*<\dots<x_m^*$ and
$\nm{x_i^*}\le 1$ for each $i$. Note that these functionals have norm
at most 1. If $k\in\N$, let $\G_k^X$ be the set of sequences
$y_1^*<\dots<y_k^*$ such that $y_i^*\in {\bf Q}$ for each $i$,
$y_1^*\in A_{j_{2k}}^*(X)$ and $y_{i+1}^*\in A_{\sigma(\seq {{y^*}}
i)}^*(X)$ for each $1\le i\le k-1$. We call these {\sl special
sequences}. Let $B_k^*(X)$ be the set of functionals of the form
$f(k)^{-1/2}\sum_{j=1}^kg_j$ such that $(\seq g k)\in\G_k^X$. These,
when $k\in K$, are the {\sl special functionals} (on $X$ of size
$k$). Note that if $g\in\cal F$ and $g(k)=f(k)^{1/2}$, then a special
functional of size $k$ is also a $(k,g)$-form.

Now, given a proper set $\cal S$ of spreads, we define the space
$X(\cal S)$, inductively. It is the completion of $c_{00}$ in the
smallest norm satisfying the following equation.
$$\eqalign{\nm x=\nm x_{c_0}&\vee
\sup\Bigl\{f(n)^{-1}\sum_{i=1}^n\nm{E_ix}:2\le n\in\N,\sleq E n\Bigr\}\cr
&\vee\sup\Bigl\{|x^*(Ex)|:k\in K,\, x^*\in B_k^*(X),\, E\subset\N\ 
\hbox{an interval}\Bigr\}\cr
&\vee\sup\{\nm{Sx}:S\in\cal S\}\cr}$$
Note that the sets $\sleq E n$ above are intervals. In the case
$\cs=\{Id_{c_{00}}\}$ the fourth term drops out and the definition
reduces to that of the space constructed in [GM].  The fourth term is
there to force $X(\cs)$ to have property (i) claimed in the
theorem. It is easy to verify that this is the case.  The second term
ensures that $X$ satisfies a lower $f$-estimate. It is also not hard
to show that $X(\cs)$ is reflexive. (A proof can be found in [GM],
end of section 3,
which works in this more general context.)

We now prove a few lemmas with the eventual aim of proving that the
spaces $X(\cs)$ have the second property claimed in the main theorem.
The first we quote from [GM].

\proclaim Lemma 6. Let $K_0\subset K$. There exists a function 
$g\in\cal F$ such that $g\ge\sqrt f$, $g(k)=\sqrt{f(k)}$ whenever
$k\in K_0$ and $g(x)=f(x)$ whenever $N\in J\setminus K_0$ and $x$ is in
the interval $[\log N,\exp N]$. \hfill $\square$

\proclaim Lemma 7. Let $0<\e\le 1$, $0\le\d<1$, $M\in L$ and let $n$ 
and $N$ be integers such that $N/n\in [\log M,\exp M]$ and
$f(N)\le(~1~+~\d~)f(N/n)$.  Assume that $\seq x N$ satisfies the
RIS($\e$) condition and let $x=\speq x N$. Then
$\bnrm{(f(N)/N)x}{n}\le(1+\d)(1+3\e)$. In particular, if $n=1$, we
have $\nm{(f(N)/N)x}\le(1+3\e)$.

\Proof  Let $g$ be the function given by Lemma 6 in the case $K_0=K$.
It is clear that every vector $Ex$ such that $\nm{Ex}>1$ is normed by
a $(k,g)$-form for some $k$, so the conditions of Lemma 3 are
satisfied. Let $\sleq E n$ be successive intervals and let $w$ be the
weight function from that proof. Then $w(\ran x)=N$ and, using
inequality $(*)$ from that proof and noting that $N/n\in[\log M,\exp
M]$, we obtain 
$$\eqalign{\sm i n\nm{E_ix}&\le(1+2\e)(N/g(N/n)+Nn/n_1)\cr
&=(1+2\e)(N/f(N/n)+Nn/n_1)\cr
&\le(1+2\e)((1+\d)N/f(N)+Nn/2^{N^2/\e^2})\cr &\le(1+3\e)(1+\d)N/f(N)\
.\cr}$$ 
This proves the lemma. \hfill $\square$ \bigskip

The key lemma used to prove that property (ii) holds is a
generalization of Lemma 22 of [GM]. Note first that a proper set $\cs$
of spreads must be countable, and if we write it as
$\{S_1,S_2,\dots\}$ and set $\cs_m=\{\seq S m\}$ for every $m$, then
for any $x\in X(\cs)$, $x^*\in X(\cs)^*$, we have
$\lim_m\,\sup\{|x^*(Ux)|:U\in\cs\setminus\cs_m\}\le
\nm x_\infty\nm{x^*}_\infty$.

\proclaim Lemma 8. Let $\cal S$ be a proper set of spreads, let 
$X=X(\cs)$, let $Y\subset X$ be an infinite-dimensional block subspace
and let $T$ be a continuous linear operator from $Y$ to $X$. Let
$\cs=\bigcup_{m=1}^\infty\cs_m$ be a decomposition of $\cal S$
satisfying the condition just mentioned. Then for every $\e>0$ there
exists $m$ such that, for every $x\in Y$ such that $\nrm xm\le 1$ and
$\supp(x)>\{m\}$,
$$d(Tx,m\conv\{\la Ux:U\in\cs_m, |\la|=1 \})\le\e\ .$$

\Proof  It is not hard to show that $T$ can be perturbed (in the 
operator norm) to an operator whose matrix (with respect to the
natural bases of $X$ and $Y$) has only finitely many non-zero entries
in each row and column. We may therefore assume that $T$ has this
property. We may also assume that $\nm T\le 1$.

Now suppose that the result is false. Then, for some $\e>0$, we can
find a sequence $\sq y$ with $y_n\in Y$, $\nrm {y_n}{n}\le 1$ and
$\supp(y_n)>\{n\}$ such that, setting $C_n=n\conv\{\la
Uy_n:U\in\cs_n, |\la|=1 \}$, we have $d(Ty_n,C_n)>\e$, and also such that if
$z_n$ is any one of $y_n$, $Ty_n$ or $Uy_n$ for some $U\in\cs_n$ and
$z_{n+1}$ is any one of $y_{n+1}$, $Ty_{n+1}$ or $Vy_{n+1}$ for some
$V\in\cs_{n+1}$, then $z_n<z_{n+1}$.

By the Hahn-Banach theorem, for every $n$ there is a norm-one
functional $y_n^*$ such that
$$\sup\{y_n^*(x):x\in C_n+\e B(X)\}<y_n^*(Ty_n)\ .$$
It follows that $y_n^*(Ty_n)>\e$ and $\sup |y_n^*(C_n)|\le 1$. 
Therefore $|y_n^*(Uy_n)|\le n^{-1}$ for every $U\in\cs_n$. We
may also assume that the support of $y_n^*$ is contained in the smallest 
interval containing the supports of $y_n$, $Ty_n$ and $Uy_n$ for
$U\in \cs_n$. (The case of complex scalars requires a standard
modification.)

Given $N\in L$ define an $N$-{\sl pair} to be a pair $(x,x^*)$
constructed as follows.  Let $y_{n_1},y_{n_2},\dots,y_{n_N}$ be a
subsequence of $\sq y$ satisfying the RIS(1) condition, which implies
that $n_1>N^2$.  Let $x=N^{-1}f(N)(y_{n_1}+\dots+ y_{n_N})$ and let
$x^*=f(N)^{-1}(y_{n_1}^*+\dots+y_{n_N}^*)$, where the $y_{n_i}^*$ are
as above. Lemma 7 implies that $\|x\| \le 4$ and $\nrm x {\sqrt N}\le 8$.

If $(x,x^*)$ is such an $N$-pair, then $x^*\in A_N^*(X)$ and,
by our earlier assumptions about supports, $$x^*(Tx)=N^{-1}\sm i
Ny_{n_i}^*(Ty_{n_i})>\e\ .$$ Similarly, $|x^*(Ux)|\le N^{-2}$ for
every $U\in\cs_N$.

Let $k\in K$ be such that $(\e/24)f(k)^{1/2}>1$. We now construct
sequences $\seq x k$ and $x_1^*,\dots,x_k^*$ as follows. Let
$N_1=j_{2k}$ and let $(x_1,x_1^*)$ be an $N_1$-pair. Let $M_2$ be such
that $|x_1^*(Ux_1)|\le\nm{x_1}_\infty\nm{x_1^*}_\infty$ if
$U\in\cs\setminus\cs_{M_2}$. The functional $x_1^*$ can be perturbed
so that it is in $\bf Q$ and so that
$\sigma(x_1^*)>\max\{M_2,f^{-1}(4)\}$, while $(x_1,x_1^*)$ is still an
$N_1$-pair. In general, after $\seq x {{i-1}}$ and
$x_1^*,\dots,x_{{i-1}}^*$ have been constructed, let $(x_i,x_i^*)$ be
an $N_i$-pair such that all of $x_i,Tx_i$ and $x_i^*$ are supported
after all of $x_{i-1},Tx_{i-1}$ and $x_{i-1}^*$, and then perturb
$x_i^*$ in such a way that, setting $N_{i+1}=\s(x_1^*,\dots,x_i^*)$,
we have $|x_i^*(Ux_i)|\le\nm{x_i}_\infty\nm{x_i^*}_\infty$ whenever
$U\in\cs\setminus\cs_{N_{i+1}}$ and we also have $f(N_{i+1})>2^{i+1}$
and $\sqrt{f(N_{i+1})}>2|\ran(\sm j ix_j)|$.

Now let $x=(\speq x k)$ and let $x^*=f(k)^{-1/2}(x_1^*+\dots+x_k^*)$.
Our construction guarantees that $x^*$ is a special functional, and
therefore of norm at most 1. We therefore have $$\nm{Tx}\ge x^*(Tx)>\e
kf(k)^{-1/2}\ .$$ Our aim is now to get an upper bound for $\nm x$ and
to deduce an arbitrarily large lower bound for $\nm T$. For this
purpose we use Lemma 3.

Let $g$ be the function given by Lemma 6 in the case
$K_0=K\setminus\{k\}$. It is clear that all vectors $Ex$ are either
normed by $(M,g)$-forms or by spreads of special functionals of length
$k$, or they have norm at most 1. In order to apply Lemma 3 with this
$g$, it is therefore enough to show that $|U^*z^*(Ex)|=|z^*(UEx)|\le
1$ for any special functional $z^*$ of length $k$ and $U\in\cal S$.
Let $z^*=f(k)^{-1/2}(z_1^*+\dots+z_k^*)$ be such a functional with
$z_j^*\in A_{m_j}^*$. Suppose that $U\in\cs_{M+1}\setminus\cs_M$, and
let $j$ be such that $N_j\le M<N_{j+1}$. Let $t$ be the largest
integer such that $m_t=N_t$. Then $z_i^*=x_i^*$ for all $i<t$,
because $\sigma$ is injective. For
such an $i$, $|z_i^*(Ux_i)|=|x_i^*(Ux_i)|<N_i^{-2}$, if $M<N_i$. If
$M\ge N_{i+1}$, then $U\notin\cs_{N_{i+1}}$, so
$|x_i^*(Ux_i)|\le\nm{x_i}_\infty\nm{x_i^*}_\infty\le 2^{-i}$. If
$N_i\le M<N_{i+1}$, the only remaining case, then $i=j$ and at least
we know that $|z_i^*(Ux_i)|\le \|x_i\| \le 4$.

If $l\ne i$ or $l=i>t$, then $z_l^*(Ux_i)=U^*z_l^*(x_i)$, and we have
$U^*z_l^*\in A_{m_l}^*$ for some $m_l$. Moreover, because $\s$ is
injective and by definition of $t$, in both cases $m_l\ne N_i$. If
$m_l<N_i$, then, as we remarked above, $\nrm{x_i}{\sqrt{N_i}}\le 8$,
so the lower bound of $j_{2k}$ for $m_1$ tells us that
$|U^*z_l^*(x_i)|\le k^{-2}$. If $m_l>N_i$, the same conclusion follows
from Lemma 1. There are at most two pairs $(i,l)$ for which $0\ne
z_l^*(UEx_i)\ne z_l^*(Ux_i)$ and for such a pair $|z_l^*(UEx_i)|\le
1$.

Putting all these facts together, we get that $|z^*(UEx)|\le 1$, as
desired. We also know that $(1/8)(\seq x k)$ satisfies the RIS(1)
condition.  Hence, by Lemma 3, $\nm x\le
24kg(k)^{-1}=24kf(k)^{-1}$.  It follows that $\nm
T\ge(\e/24)f(k)^{1/2}>1$, a contradiction.
\hfill $\square$ \bigskip

\proclaim Lemma 9. Let $\cs$, $X$, $Y$, $T$ and $\e$ be as in the 
previous lemma, let $m$ be as given by that lemma and let
$\ca_m=m\conv\{\la\cs_m : |\la|=1 \}$. Then there exists $U\in\ca_m$ such that
$\nmm{T-U}\le 17\e$.

\Proof  As in the last lemma, we can assume that the matrix of $T$ has
only finitely many non-zero entries in each row and column. If the
statement of the lemma is false, then for every $U\in\ca_m$ there is a
sequence $\ux=\ux_U\in\cl$ of vectors in $Y$ such that
$\lim\sup_n\nm{(T-U)(\ux)_n}>17\e$.  We will write this symbolically
as $\nm{(T-U)\ux}>17\e$.  Our first aim is to show that these $\ux_U$
can be chosen continuously in $U$.  (This statement will be made more
precise later.)

Let $\sqf {{\cu}} j k$ be a covering of $\ca_m$ by open sets of
diameter less than $\e$ in the operator norm. For every
$j=1,\dots,k$, let $U_j\in\cu_j$ and let $\ux_j$ be a sequence with
the above property with $U=U_j$. By the condition on the diameter of
$\cu_j$, we have $\nm{(T-U)\ux_j}>16\e$ for every $U\in\cu_j$. Let
$\sqf {\phi} j k$ be a partition of unity on $\ca_m$ with $\phi_j$
supported inside $\cu_j$ for each $j$.

Let $N\in L$ be greater than $k$ and $m^2$. For each $j\le k$, let
$x_{j,n_1},\dots,x_{j,n_N}$ satisfy the RIS(1) condition and let
$m<x_{j,n_1}$. Let $y_j=N^{-1}f(N)(x_{j,n_1}+\dots+x_{j,n_N})$. Let
this be done in such a way that $\sleq y k$ and also
$(T-U)x_{j,n_1}<\dots<(T-U)x_{j,n_N}$ for every $j$ and every
$U\in\ca_m$.  Finally, let the $x_{j,n_i}$ be chosen so that
$\nm{(T-U)x_{j,n_i}}>16\e$ for every $U\in\cu_j$.

Now let us consider the vector $y(U)=\sm j k\phi_j(U)y_j$.  By Lemma 7
we know that, for each $y_j$, $\nrm{y_j}{\sqrt N}\le 8$, from which it
follows that $\nrm{y(U)}{\sqrt N}\le 8$. We shall show that $y(U)$ is
a ``bad'' vector for $U$, by showing that $\nm{(T-U)y(U)}>8\e$.

To do this, let $U\in\ca_m$ be fixed and let $J=\{j:\phi_j(U)>0\}$.
Note that $\nm{(T-U)\ux_j}>16\e$ for every $j\in J$. For such a $j$
and for $i\le N$ let $z_{j,i}^*$ be a norm-one functional such that
$z_{j,i}^*\bigl((T-U)x_{j,n_i}\Bigr)>16\e$. Let these functionals be
chosen to be successive. Let
$z_j^*=f(N)^{-1}(z_{j,1}^*+\dots+z_{j,N}^*)$ and $z^*=\sum_{j\in
J}z_j^*$. Then $z_j^*(T-U)y_j>16\e$, so
$$z^*\Bigl((T-U)y(U)\Bigr)=z^*\Bigl(\sum_{j\in J}
\phi_j(U)(T-U)y_j\Bigr)>16\e .$$
However, $\nm{z^*}\le f(kN)/f(N)\le 2$, proving our claim.

The function $U\mapsto y(U)$ is clearly continuous. The
vector $y(U)$ satisfies $m<y(U)$, $\nrm{y(U)}m\le 8$, and
$\nm{(T-U)y(U)}>8\e$. We now apply a fixed-point theorem.

For every $U\in\ca_m$, let $\Gamma(U)$ be the set of $V\in\ca_m$ such
that $\nm{(T-V)y(U)}\le 8\e$. Clearly $\Gamma(U)$ is a compact convex
subset of $\ca_m$. By the previous lemma, $\Gamma(U)$ is non-empty for
every $U$. The continuity of $U\mapsto y(U)$ gives that $\Gamma$ is
upper semi-continuous, so there exists a point $U\in\ca_m$ such that
$U\in\Gamma(U)$. But this is a contradiction.
\hfill $\square$ \bigskip

Lemma 9 shows in particular that any operator $T:Y\ra X$ can be 
approximated arbitrarily well in the $\nmm.$-norm by the restriction
of some operator $U\in\ca$. We have therefore finished the proof of
property (ii). The proof of (iii) is much easier, and will complete
the proof of Theorem 5.

\proclaim Lemma 10. The seminorm $\nmm.$ satisfies the algebra
inequality $\nmm{UV}\le\nmm U\nmm V$.

\Proof  This lemma is one of the main reasons for the technical
differences between this paper and [GM]. To see it, pick $\e>0$ and
let $\sq x\in\cl$ be a sequence such that
$\nm{UVx_n}\ge(1+\e)^{-1}\nmm{UV}$ for every $n$. After suitable
perturbations and selections of subsequences we may assume that the
sequences $\sq x$, $\sq {Vx}$ and $\sq {UVx}$ are successive.

Given $m\in L$ sufficiently large, construct a vector $(y_m)$ as
follows. Choose $k\ge 8/\e$, let $M=m^k$ and let
$x_{n_1},\dots,x_{n_M}$ satisfy the $RIS(\e)$ condition. Then let
$y_m={f(M)\over M}(x_{n_1}+\dots+x_{n_M})$. For $0\le i<m^2$ and
$0\le j<m^2$, we now let
$$z_{ij}={f(M)\over M}
\sum_{s=im^{k-2}+jm^{k-4}+1}^{im^{k-2}+(j+1)m^{k-4}}x_{n_s}\qquad\qquad
\hbox{and}\qquad\qquad z_i=\sum_{j=0}^{m^2-1}z_{ij}\ .$$ 
By Lemma 7 we know that 
$$\nm{z_{ij}}\le(1+3\e){m^{k-4}\over
f(m^{k-4})}{f(M)\over M}\le(1+\e)(1+3\e)m^{-4}$$ for each $i$ and
$j$. It follows from the proof of Lemma 4 that $\nrm {z_i} {m}\le
m^{-2}(1+\e)(1+3\e)(1+2m/m^2)$ for each $i$. If $m$ is large enough,
we then have $\nm{Vz_i}\le m^{-2}(1+\e)^2(1+3\e)\nmm V$ for every $i$.
By the proof of Lemma 4 again, we find that $\nrm {Vy_m}{m}\le
(1+\e)^2(1+3\e)(1+2m/m^2)\nmm V$, and hence, for $m$ large enough,
$\nm{UVy_m}\le(1+\e)^3(1+3\e)\nmm U\nmm V$.

On the other hand, for every $m$
$$\nm{UVy_m}={f(M)\over M}\Bnm{\sm i MUVx_{n_i}}
\ge{1\over M}\sm i M\nm{UVx_{n_i}}\ge(1+\e)^{-1}\nmm{UV}\ .$$
Letting $\e$ tend to zero we obtain the desired inequality. 
\hfill $\square$ \bigskip

\noindent {\bf \S 4. Applications}
\medskip

\noindent {\bf (4.1)}\ Let $\cs=\{Id\}$, let $X=X(\cs)$, let $Y$ be
any block subspace of $X$ and let $i$ be the inclusion map from $Y$ to
$X$. Then given any operator $T$ from $Y$ to $X$, there exists by
Theorem 5, for every $\e>0$, some $\la$ such that $\nmm{T-\la
i}<\e$. Since $|\la|\le\nmm T+\e$, an easy compactness argument then
shows that there exists $\la$ such that $\nmm{T-\la i}=0$ and thus that
$T-\la i$ is strictly singular, which is one of the main results of
[GM]. It implies easily that $X$ contains no unconditional basic
sequence, is hereditarily indecomposable and is not isomorphic to any
proper subspace. (The third fact is true because every operator on $X$
must be strictly singular or Fredholm with index zero.) In this case
the algebra $\cal G$ (defined just after Theorem 5) reduces to the
field of scalars $\R$ or $\C$.  It follows that $X^n$ is isomorphic to
$X^m$ if and only if $m=n$. Indeed, if $n>m$, the image under $\phi$
of any $T\in L(X^n,X^m)$ is a rectangular matrix $A\in M_{m,n}$ which
has non-zero kernel. It follows easily that $T$ is singular. \bigskip

\noindent {\bf (4.2)}\ Let $\cal S$ be the proper set mentioned
earlier, generated by the shift, which we denote by $S$. That is,
$\cal S$ consists of all maps of the form $S_{A,B}$ where
$A=[m,\infty)$ and $B=[n,\infty)$. We will write $L$ for the left
shift, which is (formally) the adjoint of $S$. Then every operator in
$\cal S$ is of the form $S^mL^n$. Since $SL-I$ is of rank one, every
operator in $\cal A$ is a finite-rank perturbation of an operator of
the form $\sum_{n=0}^N\la_nS^n+\sm n N\mu_nL^n$, so the difference is
of $\nmm.$-norm zero.

\proclaim Lemma 11. Let $U=\sum_{n=0}^N\la_nS^n+\sm n N\mu_nL^n$. Then 
$\nm U=\nmm U=\sum_{n=0}^N|\la_n|+\sm n N|\mu_n|$.

\Proof  Note that $\sum_{n=0}^N|\la_n|+\sm n N|\mu_n|$ is the norm of 
$U$ considered as an operator on $\ell_1$. Clearly it is enough to
prove the inequality $\nmm U\ge\nm U_{\ell_1\ra\ell_1}$. For
notational convenience, let $\la_{-n}=\mu_n$ for $1\le n\le N$.

For an integer $r\in L$ consider the vector $x_r=\sm j r
\bfe_{3jN}$. Since every unit vector $\bfe_i$ satisfies 
$\nrm{\bfe_i}{n}=1$ for
every $n$, we have $\nm{x_r}\le r/f(r)$, by Lemma 7. On the other
hand, splitting $Ux_r$ into $3rN$ singleton pieces from $N+1$ to
$(3r+1)N$ gives that $\nm {Ux_r}\ge (r/f(3Nr))\sum_{n=-N}^N|\la_n|$. As
$r\ra\infty$, $f(r)/f(3Nr)\ra 1$, which shows that $\nm U$ satisfies
the required inequality. To get it for $\nmm U$, let $\e>0$ and
$n\in\N$. Then, by Lemma 7, there exists $r$ such that $\nrm
{(f(r)/r)x_r}{n}\le 1+\e$.  This is also true if $x_r$ is shifted, so
the lower bound above on $\nm{Ux_r}$ gives the result. \hfill
$\square$ \bigskip

Since all powers of $S$ and $L$ have norm 1, it is obviously true that
every operator of the form $\sum_{n=0}^\infty a_nS^n+\sum_{n=1}^\infty
a_{-n}L^n$ is continuous if $\sum_{n=-\infty}^\infty
|a_n|<\infty$. The next result gives, up to a strictly singular
perturbation, the converse of this fact.

\proclaim Corollary 12. There is an algebra homomorphism and
projection $\phi:L(X)\ra L(X)$ onto the subspace consisting of
Toeplitz operators with absolutely summable coefficients. If
$T\in L(X)$ then $\nmm{\phi(T)-T}=0$.

\Proof Recall the remark following the statement of Theorem 5. In this
case, by Lemma 11, the algebra $\cal G$, the $\nmm.$-completion of
$\ca$, is the same as the completion in $L(X)$ and also the completion
in the operator norm on $\ell_1$. Therefore $\cal G$ can be regarded
as a subalgebra of $L(X)$ consisting of Toeplitz operators with
absolutely summable coefficients. If we do this, then the algebra
homomorphism $\phi$ defined after Theorem 5 is also a projection.
The equation $\nmm{\phi(T)-T}=0$ follows easily from the definition
of $\phi$. \hfill $\square$ \bigskip

Recall that a Banach space is said to be {\sl prime} if it is
isomorphic to every infinite-dimensional complemented subspace of
itself. The only known examples before this paper were $c_0$ and
$\ell_p$ ($1\le p\le\infty$). These were shown to be prime by
Pe\l czy\'nski [P], apart from $\ell_\infty$ which is due to Lindenstrauss
[L].  The space $X$ is prime by virtue of having no non-trivial
complemented subspaces and being isomorphic to its subspaces of finite
codimension.

\proclaim Theorem 13. The space $X$ is prime.

\Proof Let $P:X\ra X$ be a projection. By the previous corollary the
operator $\phi(P)$ is a convolution by some absolutely summable
sequence $(a_n)_{n\in\Z}$. Moreover, $\phi(P)^2=\phi(P)$. But the
Fourier transform of the sequence $(a_n)_{n\in\Z}$ is a continuous
function on the circle squaring to itself. Hence it is constantly zero
or one. It follows that $a_0$ is zero or one and all the other $a_n$
are zero. That is, $\phi(P)$ is zero or the identity.  Since
$P-\phi(P)$ is strictly singular, it follows that $P$ is of finite
rank or corank. Thus, if $PX$ is infinite-dimensional, then it has
finite codimension. Since the shift on $X$ is an isometry, it follows
that $X$ and $PX$ are isomorphic, which proves the theorem. \hfill
$\square$ \bigskip

We note here that the argument in the above proof can be generalized
to show that if $m$ and $n$ are integers with $m>n$, then $X^n$ does
not contain a family $\seq P m$ of infinite-rank projections
satisfying $P_iP_j=0$ whenever $i\ne j$. Indeed, given any projection
$P\in L(X^n)$, we can regard it as an element of $M_n(L(X))$. Acting
on each entry with first $\phi$ and then the Fourier transform, we get
a function $h\in M_n(C(\T))$. The map taking $P$ to $h$ is an algebra
homomorphism so $h$ is an idempotent. Regarding $h$ as a continuous
function from $\T$ to $M_n(\C)$, we have that $h(t)$ is an idempotent
in $M_n(\C)$ for every $t\in\T$. By the continuity of rank for
idempotents, we have that if $h(t)=0$ for some $t$, then $h$ is
identically zero. But then $P$ is strictly singular and hence of
finite rank. Applying this reasoning to the family $\seq P m$ above,
we obtain $\seq h m$ such that, for every $t\in\T$,
$h_1(t),\dots,h_m(t)$ is a set of non-zero idempotents in $M_n(\C)$ with
$h_i(t)h_j(t)=0$ when $i\ne j$. But this is impossible if $m>n$.
It follows that $X^n$ and $X^m$ are isomorphic if and only if
$n=m$. 

Another simple consequence of Corollary 12 is that, up to strictly
singular perturbations, any two operators on $X$ commute. Indeed, if
$V$ and $W$ are two operators, then $\phi(V)$ and $\phi(W)$ commute,
so $\phi(VW-WV)=0$, from which it follows that $\nmm{VW-WV}=0$.

For the rest of this section, we assume that $X$ has complex scalars.
Let $\psi:L(X)\ra C(\T)$ be the composition of $\phi$ with the Fourier
transform. Then $\psi$ is also a continuous algebra homomorphism.
Given an operator $T$, let $K_T$ be the compact set of $\mu\in\C$ such
that $\mu$ is infinitely singular for $T$. (Recall that this means
that for every $\e>0$ there is an infinite-dimensional subspace
$Y\subset X$ such that $\nm{Ty-\mu y}\le\e\nm y$ for every $y\in Y$.)
Since $T-\phi(T)$ is strictly singular, $K_{\phi(T)}=K_T$.

\proclaim Lemma 14. The function $\psi(T)$ takes the value zero at
some $\exp(i\theta)$ if and only if $0$ is infinitely singular for
$T$.

\Proof If $\psi(T)$ takes the value zero at $\exp i\theta$, we can
construct an approximate eigenvector for $\phi(T)$ with eigenvalue
zero as follows. Suppose that $\phi(T)$ is convolution by the sequence
$(a_n)_{n\in\Z}$, and let $\e>0$. We know that
$\sum_{n\in\Z}a_n\exp(in\theta)=0$. Let $N\in L$ and let $x_N$ be the
vector $(f(N^2)/N^2)\sum_{n=N^2}^{2N^2}\exp(in\theta) \bfe_n$. By Lemma 7
we have $\nm{x_N}=1$. Let $U$ be convolution by the sequence
$(a_n)_{n=-N}^N$. If $N$ is large enough, then
$\nm{U-\phi(T)}\le\e/2$, since $(a_n)_{n\in\Z}$ is absolutely
summable.  Moreover, all but at most $4N$ of the possible $N^2+2N$
non-zero coordinates of $Ux_N$ are equal to
$(f(N^2)/N^2)\sum_{n=-N}^Na_n\exp(in\theta)$. Taking $N$ sufficiently
large, we can therefore make $\nm{\phi(T)-U}$ and $\nm{Ux_N}$ as small
as we like. Therefore zero is infinitely singular for $\phi(T)$. Since
$\nmm{T-\phi(T)}=0$, the same is true for~$T$.

Conversely, if $\psi(T)$ never takes the value zero, then it can be
inverted in $C(\T)$. A classical result states that the Fourier
transform of this inverse will also be in $\ell_1(\Z)$, so in
particular $\phi(T)$ has an inverse $U$ which is continuous when
considered as an operator on $X$ and satisfies $U=\phi(U)$. Therefore
$\phi(UT-I)=0$, so $UT-I$ is strictly singular and zero is not
infinitely singular for $T$. \hfill $\square$ \bigskip

\proclaim Corollary 15. $K_T$ is the image under $\psi(T)$ of the
unit circle $\T$. 

\Proof  This follows from Lemma 14 applied to the operator $T-\la$.
\hfill $\square$ \bigskip

\proclaim Theorem 16. A subspace $Y$ of $X$ is isomorphic to
$X$ if and only if it has finite codimension.

\Proof Let $T:X\ra Y$ be an isomorphism. Then 0 is not infinitely
singular for $T$, so, as in the proof of Lemma 14, we can find $U$
such that $TU$, $UT$ and $I$ are the same, up to a strictly singular
perturbation. Since $TU-I$ is strictly singular, $TU$
is Fredholm with index zero. In particular $\codim Y=\codim
TX\le\codim TUX<\infty$. As we have already mentioned, the if part
follows from the existence of the isometric shift. \hfill $\square$ \bigskip

\bigskip

\noindent {\bf (4.3)} Let $\cs$ be the proper set generated by the
double shift $S^2$. That is, $\cs$ is as in the previous example but
$m$ and $n$ are required to be even. We adapt a result about Fredholm
operators to show that every operator on $X(\cs)$ has even
index. Suppose that this is not true and let $T$ be an operator of odd
index. By Theorem 5, and by the fact that every operator in $\cs$
differs by a finite-rank operator from some even shift, we can find,
for any $\e>0$, some linear combination $U$ of even shifts such that
$\nmm{T-U}<\e$. We obtain a contradiction by showing that no such $U$
can have odd index and that there is an $\nmm.$-neighbourhood of $T$
inside which all operators do have odd index.

\proclaim Lemma 17. Let $U$ be a Fredholm isometry on a Banach space 
$X$ with a left inverse $V$, and let $T:X\ra X$ be a Fredholm operator
which can be written in the form $P(U)+Q(V)$ for polynomials $P$ and
$Q$. Then the index of $T$ is a multiple of the index of $U$.

\Proof Assume first that the scalars are complex.
Given any operator $W$ on $X$, define $F_W$ as in [GM] to be
the set of $\la\in\C$ such that $W-\la$ is an isomorphism on some
finite-codimensional subspace. (It is not hard to show that $F_W$ is
the complement of $K_W$ defined in the last section.) Then $F_W$ is
open and $(W-\la)(X)$ is closed,
$\dim\ker(W-\la)<\infty$ whenever $\la\in F_W$. Hence, the
operator $W-\la$ is quasi-Fredholm, and the generalized index is
constant on connected components of $F_W$.

In the case of the operator $U$, it is clear that $F_U$ contains the
open unit disc and the set of all $\la$ such that $|\la|>1$. Hence,
either $F_U$ is connected and the index of $U-\la$ is constantly zero
on $F_U$, or $\C\setminus F_U=\T$. In the second case, the index of
$U-\la$ is zero if $|\la|>1$ and $\ind U$ if $|\la|<1$. In either case,
the only possible values are $0$ and $\ind(U)$.

Now suppose that $T$ is as in the statement of the lemma. For
sufficiently large $N$, $TU^N$ can be written $R(U)$ for some
polynomial $R$ and is still Fredholm. Writing $R(U)=c\prod_i(U-\la_i)$,
we must have $\la_i\in F_U$ for $TU^N$ to be Fredholm, so $\ind(U-\la_i)$
is either $0$ or $\ind(U)$. It follows that the index of $R(U)$, and 
hence that of $T$, is a multiple of the index of $U$ as stated.

 The real case follows by considering the extension $U_\C$ of $U$
as an isometry on the complexification $X_\C = X\oplus_2 X$ of $X$.
\hfill $\square$ \bigskip 

We now use the following lemma, which is an easy variant of a standard
lemma which can be found, for example, as [LT Prop. 2.c.9].

\proclaim Lemma 18. Let $X$ and $Y$ be any Banach spaces and let 
$T:X\ra Y$ be a Fredholm operator. There exists $\e>0$ such that if
$S$ is any operator such that every infinite-dimensional
subspace of $X$ contains some $x$
for which $\nm{Sx}<\e\nm x$, then $T+S$ is Fredholm of the same index
as $T$.

\Proof  Pick $X_1$ such that $X=X_1\oplus\ker T$, so $T_1=T|_{X_1}$ 
is an isomorphism and let $\e<(1/2)\nm{T_1^{-1}}^{-1}$. If $T+S$ has
infinite-dimensional kernel, then so does $T_1+S_1$ (where
$S_1=S|_{X_1}$), so there is an infinite-dimensional $Z\subset X_1$ on
which $T_1+S_1=0$. This contradicts our choice of $\e$. If $(T+S)$ 
does not have closed range, then it is standard that there exists
an infinite-dimensional subspace $Z\subset X_1$ on which
$\nm{(T+S)|_Z}<\e$. This again contradicts our choice of $\e$.

Hence $\ind(T+tS)$ is defined for all $t\in[0,1]$. This function is
known to be continuous and therefore $\ind(T)=\ind(T+S)$.
\hfill $\square$ \bigskip

Putting these facts together, we find that no continuous operator on
$X(\cs)$ can be Fredholm with odd index. We therefore have the
following result.

\proclaim Theorem 19. The space $X(\cs)$ is isomorphic to its subspaces
of even codimension while not being isomorphic to those of odd
codimension. In particular, it is isomorphic to its subspaces of
codimension two but not to its hyperplanes.

\noindent Remarks: A slight modification of the above approach is also
possible. As in the last section, one can obtain a continuous algebra
homomorphism from $L(X)$ to the Toeplitz operators on $X$
corresponding to sequences in $\ell_1(\Z)$. The proof of Theorem 19 is
then as above except that instead of Lemma 18 one can use the lemma of
which it is a variant, which states that a sufficiently small
perturbation in the operator norm does not change the index of a
Fredholm operator. The proof of Theorem 13 gives for this space also
that every complemented subspace has finite dimension or codimension.
Combining this observation with Theorem 19, we see that the space has
exactly two infinite-dimensional complemented subspaces, up to
isomorphism. It is true for this space as well that it is isomorphic
to no subspace of infinite codimension.

Note that the methods of this section generalize easily to proper sets
generated by larger powers of the shift. For example, there is a space
$X$ such that two finite-codimensional subspaces are isomorphic if and
only if their codimensions are equal~mod~7.

\bigskip

\noindent {\bf (4.4)}\ This application is more complicated than the
previous ones. The aim is to construct a space $X$ which is isomorphic
to $X\oplus X\oplus X$ but not to $X\oplus X$. There is a very natural
choice of $\cs$ in this case. For $i=0,1,2$ let $A_i$ be the set of
positive integers equal to $i+1$ ($\mod 3$), let $S_i'$ be the spread
from $\N$ to $A_i$ and let $\cs'$ be the semigroup generated by
$S_0'$, $S_1'$ and $S_2'$ and their adjoints. We shall show later that
this is a proper set. The space $X(\cs')$ is easily seen to be
isomorphic to its cube, and this isomorphism is achieved in a
``minimal'' way. (The primes in this paragraph are to avoid confusion
later.)

We shall indeed consider the space $X(\cs')$ defined above. However,
we define it slightly less directly, which helps with the proof later
that it is not isomorphic to its square. The algebra $\ca'$ arising
from the above definition is, if completed in the $\ell_2$-norm,
isometric to the Cuntz algebra $\co_3$, which was analysed using
K-theory by Cuntz in [C]. Our proof is inspired by his paper, although
K-theory is not mentioned explicitly. We shall sketch a more directly
K-theoretic approach at the end of the section.

Let $\cal T$ be the ternary tree $\bigcup_{n=0}^\infty \{0,1,2\}^n$.
Let $Y_{00}$ be the vector space of finitely supported scalar
sequences indexed by $\cal T$ (including the empty sequence).  Denote
the canonical basis for $Y_{00}$ by $(e_t)_{t\in \cal T}$, write $e$
for $e_\emptyset$, and denote the length of a word $t\in \cal T$ by $|t|$.
If $s,t\in \cal T$, let $(s,t)$ stand for the concatenation of $s$
and $t$.  We shall now describe some operators on $Y_{00}$.  

Let $S_i$ and $T_i$, for $i=0,1,2$ be defined by their action on the 
basis as follows:
$$ S_i e_t = e_{(t,i)}, \ T_i e_t = e_{(i,t)}.$$

Thus $S_i$ can be thought of as the map taking each vertex of $\ct$ to
the $i^{\rm th}$ vertex immediately below it, while $T_i$ takes the
whole tree on to the $i^{\rm th}$ branch.  The adjoints $S^*_i$ and
$T^*_i$ act in the following way: $S^*_i e_t = e_s$ if $t$ is of the
form $t=(s,i)$, and $S^*_i e_t =0$ otherwise, while $T^*_i e_t = e_s$
if $t = (i,s)$, and $T^*_i e_t = 0$ otherwise. The following facts are
easy to check: $S_i T_j = T_j S_i$, $S^*_i S_j = T^*_i T_j =
\delta_{i,j}I$; $S_i S^*_i$ and $T_i T^*_i$ are projections; if $P$
denotes the natural rank one projection on the line $\C e$, then
$\sum_{i=0}^2 S_i S^*_i = \sum_{i=0}^2 T_i T^*_i = I - P$.  Let $\cs$
and $\ca$ be respectively the proper set generated by $S_0$, $S_1$ and
$S_2$, and the algebra generated by this proper set. (Strictly
speaking, $\cs$ is not a proper set, but it is easy to embed $\ct$
into $\N$ so that the maps $S_0$, $S_1$ and $S_2$ become spreads as
defined earlier. One can use the above relations to check the
technical condition, but we do not need this. Note that $\cs$ is
the semigroup generated by the $S_i$ and the $S_i^*$, that it contains
$I$ and that $\ca$ contains $P$, as we have just shown.)
\medskip

In order to obtain the space $X$, consider the subset $\ct'$ of $\ct$
consisting of all words $t\in\cal T$ that do not start with $0$
(including the empty sequence). We modify the definition of $S_0$
slightly, by letting $S_0' e$ equal $e$ instead of $e_0$. Operators
$S_1'$ and $S_2'$ are defined exactly as $S_1$ and $S_2$ were. We
still have that the $S_i'{S_i'}^*$ are projections and that
${S_i'}^*S_j'=\d_{i,j}I$, but this time
$\sum_{i=0}^2S_i'{S_i'}^*=I$. To each $s =(i_1,\ldots,i_n) \in\ct'$ we
can associate the integer $n_s = 3^{n-1} i_1 + \cdots + 3 i_{n-1} +
i_n+1$, (with $n_\emptyset=1$), and this defines a bijection between
$\ct'$ and $\N$. The  operators $S_0'$, $S_1'$ and $S_2'$ then coincide
with the spreads on $c_{00}$ defined earlier (in fact,
$S_i' \bfe_n= \bfe_{3n-2+i}$), so we can define $\cs'$ to be the proper set
they generate and obtain the space $X(\cs')$. Let $\ca'$ be the
algebra generated by $\cs'$. We now check that $\cs'$ is a proper set
by verifying the technical condition from the definition. Using the
relation ${S_i'}^*S_j'=\d_{i,j}I$ one finds that every element of
$\cs'$ can be written in the form $UV$ where $U=S_{i_1}'\dots
S_{i_k}'$ and $V={S_{j_1}'}^*\dots {S_{j_l}'}^*$. Fix integers $m<n$.
If $l$ is larger than $\log_3(n-m)$ then at least one of $V \bfe_m$ and
$V \bfe_n$ is zero. Moreover, if $UV \bfe_m=\bfe_r$ and $UV\bfe_n=\bfe_s$, then
$U^* \bfe_r=V\bfe_m$, $U^* \bfe_s=V\bfe_n$ and for the same reason one of them is
zero if $0 < |s-r|<3^k$. The technical condition follows easily.
\medskip

For $t\in\ct$, we define $S_t$ inductively by $S_{(t,i)} = S_i
S_t$. (We also let $I=S_\emptyset$.) Thus $S_t$ takes a vertex of the
tree and moves it down the path from that vertex corresponding to the
word $t$, or in other words, $S_t e_u = e_{(u,t)}$. Let $S^*_t$ be the
adjoint of $S_t$. Then every $U\in\ca$ has a decomposition
$$U=\sum_{l=1}^N a_lS_{\alpha_l} S^*_{\beta_l}\ ,$$ 
where $\alpha_l$ and $\beta_l$ are words in $\cal T$. Define $c(U)$ to
be the smallest value of $\max_l |\beta_l|$ over all such
representations of $U$.  If $|t|>c(U)$, $U$ is decomposed as above
with $\max_l|\b_l|<|t|$  and $G_t$ is the set
of $l$ such that $t=(\g_l,\b_l)$ for some $\g_l$, then
$$ U e_t = \sum_{l\in G_t} a_le_{(\g_l,\a_l)}.$$ 
(Observe that $|\gamma_l| >0$ since $|t|>c(U)$.)

We make the obvious modifications to the above definitions for
$\ca'$. The remarks are still valid, except that the actions of
$S_t$ and $S_t'$ on $e$ will be different if the word $t$ begins
with 0. The next lemma is similar to Lemma 11.

\proclaim Lemma 20. Let $U\in\ca'$. Then for $|t|> c(U)$, we have
the inequality $\|U e_t\|_1\le\nmm U$.

\Proof Let $|t|>c(U)$ and suppose that $Ue_t = \sum_{k=1}^M c_k
e_{s_k}$, where the $s_k$s are distinct.  Since $|t|>c(U)$, we have $U
e_{u,t} = \sum_{k=1}^M c_k e_{u,s_k}$ for every $u$. Pick a sequence
$(u_j)_{j=1}^\infty$ lacunary enough to guarantee that the sequences
$(e_{u_j,t})_{j=1}^\infty$ and $(U e_{u_j,t})_{j=1}^\infty$ are
successive. Then by the construction of $X$, we obtain the inequality
$$ \| U(\sum_{j=1}^N e_{u_j,t})\| \ge {N\over f(MN)} \sum_{k=1}^M
|c_k|,$$
while  for $N\in L$
$$ \|\sum_{j=1}^N e_{u_j,t}\| \le {N\over f(N)}$$
by Lemma 7. Letting $N\rightarrow \infty$, this gives $\sum_{k=1}^M
|c_k| \le \|U\|$.  To get the inequality for $\nmm U$, we also use Lemma
7. Given $n$ and $\e>0$, it guarantees the existence of $N\in L$ such
that $\nrm {\sum_{j=1}^Ne_{u_j,t}}{n}\le(1+\e)(N/f(N))$, which is
enough. \hfill $\square$ \bigskip

We now consider the algebra $\ca$.  Let $Y$ be the completion of
$Y_{00}$ equipped with the $\ell_1$ norm (or in other words let
$Y=\ell_1(\ct)$) and let $\ce$ denote the norm closure of $\ca$ in
$L(Y)$. Note that every $S_i$ or $T_i$ is an isometry on $Y$.

\proclaim Lemma 21. Every Fredholm operator in $\ce$ has index 0. More
generally, every Fredholm operator $T:Y^q\ra Y^p$ given by a matrix in
$M_{p,q}(\ce)$ satisfies $2\ind(T)=p-q$. In particular, no operator
in $M_{p,q}(\ce)$ is Fredholm if $p-q$ is odd.

\Proof  Since the Fredholm index is stable under small perturbations, 
it is enough to consider operators in $\ca$ (as operators on $Y$). For 
any such operator $U$ we associate the operator
$$\tilde U=\sum_{i=0}^2T_iUT^*_i.$$
We claim that $\tilde U$ is a finite rank perturbation of $U$. It is enough
to show that $\tilde S_i$ is a rank-one perturbation of $S_i$ (and that
$(U^*)^{\tilde {}}=(\tilde U)^*$). But
$$\tilde S_j=\sum_{i=0}^2T_iS_jT^*_i=S_j\sum_{i=0}^2T_iT^*_i=
S_j(I-P)=S_j-S_jP\ .$$
Consider the projections $Q_i = T_i T^*_i$. Then $Q_i Q_j = 0$ 
for $i\ne j$ and
$$ Y = \C e \oplus Q_0 Y \oplus Q_1 Y \oplus Q_2 Y.$$
Each $T_i U T^*_i$ represents an operator on $Q_i Y$, equivalent (in
the obvious sense) to $U$ on $Y$, and $\tilde U$ is $0$ on the
component $\C e$. It follows that $\ind(\tilde U) = 3\ind(U)$.  On the
other hand $\ind(\tilde U) =\ind(U)$ since it is a finite rank perturbation
of $U$. It follows that $\ind(U) = 0$.

The proof is essentially the same for the more general
statement. Given a matrix $A\in M_{p,q}(\ce)$, use the tilde operation
on each entry. The resulting matrix is equivalent to three copies
of $A$ plus the zero matrix in $M_{p,q}$. This zero matrix contributes
$q$ to the dimension of the kernel and $p$ to the codimension of the
image, from which we obtain the equation
$$\ind(T)=\ind(\tilde T)=3\ind(T)+q-p\eqno\square$$ 
\bigskip

Let $\ci$ be the closed two-sided ideal in $\ce$ generated by
$P$. This ideal contains all rank-one operators of the form
$e_s^*\otimes e_t$ with $s,t\in\ct$. Hence, every finite rank operator
on $Y^n$ which is $w^*$ continuous (considering $Y^n$ as the dual of
$(c_0)^n$) belongs to $M_n(\ci)$. Indeed, the matrix of such an
operator consists of entries which are finite sums of the form $\sum_k
y_k\otimes x_k$, with $y_k\in c_0$.  We can approximate $y_k$ and
$x_k$ by finitely supported sequences $\tilde y_k$ and $\tilde x_k$,
and $\sum_k \tilde y_k \otimes \tilde x_k$ certainly belongs to
$\ci$. (In fact, $\ci$ consists exactly of the compact
$w^*$-continuous operators on $\ell_1$.)

\proclaim Lemma 22. If $V\in M_n({\cal E})$ is Fredholm then there
exists $W\in M_n(\ci)$ such that $V+W$ is invertible in $M_n(\ce)$.  

\Proof By Lemma 21 the index of $V$ is zero. Let $x_1,\ldots,x_N$ and
$z_1,\ldots,z_N$ be bases for the kernel and cokernel. We can
construct a $w^*$-continuous projection $\sum_{k=1}^N y_k \otimes x_k$
on the kernel.  Then $W = \sum_{k=1}^N y_k \otimes z_k$ will do.
\hfill $\square$ \bigskip

Let $\cal O$ denote the quotient algebra ${\cal E}/\ci$.  

\proclaim Lemma 23. Any lifting in $M_n(\ce)$ of an invertible element
in $M_n({\cal O})$ is Fredholm on $Y^n$.

\Proof  Let $xy=yx=1$ in $M_n({\cal O})$ and let $u,v$ be
any liftings of $x$ and $y$. Then $uv$ and $vu$ are compact
perturbations of the identity and hence Fredholm. It follows that
$u$ and $v$ are isomorphisms on finite codimensional subspaces and have
finite dimensional cokernels. Hence, $u$ is Fredholm.
\hfill $\square$ \bigskip

As an immediate consequence of the preceding two lemmas we have the
following statement.

\proclaim Corollary 24. Every invertible element of $M_n({\cal O})$
can be lifted to an invertible element of $M_n({\cal E})$.
\hfill $\square$

Now recall the remarks following Theorem 5. It follows easily from
Lemma 20 that $\nmm.$ is actually a norm on $\ca'$, so the Banach
algebra $\cal G$ is the $\nmm.$-completion of $\ca'$. Recall that
there is a unital algebra homomorphism $\phi:L(X)\ra\cal G$.

\proclaim Lemma 25. There is a norm-one algebra homomorphism $\theta$
from $\cal G$ to $\co$.

\Proof Define a map $\theta_0:\ca'\ra\co$ as follows. Given
$U\in\ca'$, write $U=\sm l Na_lS_{\a_l}'{S_{\b_l}'}^*$ in some way,
consider the corresponding sum $\sm l Na_lS_{\a_l}S_{\b_l}^*$ as an
element of $\ce$ and let $\theta_0(U)$ be the image of this operator
under the quotient map from $\ce$ to $\co$. To see that this map is
well defined, observe that for any pair of words $\a$ and $\b$ we have
the equation $S_\a'{S_\b'}^*=\sum_{i=0}^2S_{(i,\a)}'{S_{(i,\b)}'}^*$. If
$n$ is sufficiently large, we can therefore write $U$ as above in such
a way that all the $\a_l$ are words of length $n$. Let $W_n$ be the
set of all words of length $n$. Then what we have said implies that $U$
can be written as a sum $\sum_{\a\in W_n}S_\a'V_\a^*$, where each
$V_\a^*$ is some linear combination of distinct operators of the form
${S_{\b}'}^*$.  It is easy to see now that $U=0$ if and only if $V_\a=0$
for every $\a\in W_n$, and moreover that distinct ${S_\b'}^*$ are
linearly independent. Therefore any $U\in\ca'$ has at most one
representation in the above form. In $\ca$ we know that for any pair
of words $\a$ and $\b$ the images in $\co$ of the operators
$S_\a S_\b^*$ and $\sum_{i=0}^2S_{(i,\a)}S_{(i,\b)}^*$ are the same. It
follows that $\theta_0$ is well defined. Similarly, one can show that
it is a unital algebra homomorphism.

Let $P_n$ denote the projection on to the first $n$ levels of the tree
$\ct$, so that $P_n\in\ci$ for every $n$.  If $U\in\ca'$, then Lemma
20 implies that 
$$\lim_n \|U(I-P_n)\|_{L(Y)} \le \nmm U\ .$$ 
It follows that we may extend $\theta_0$ to a norm-one homomorphism
$\theta:\cal G\ra\co$, as claimed. \hfill $\square$ \bigskip

\proclaim Theorem 26. The spaces $X$ and $X\oplus X$ are not isomorphic.

\Proof Suppose that they were. We would then be able to find $U\in
M_{2,1} (L(X))$, $U:X\ra X\oplus X$ with an inverse $V\in
M_{1,2}(L(X))$, $V:X\oplus X\ra X$. The matrix $\left(\matrix{U&0\cr
0&V}\right)$ is an invertible element of $M_3(L(X))$ and therefore has
an invertible image in $M_3(\co)$ under $\theta\circ\phi$. By
Corollary 24 we can lift this image to an invertible element
$\left(\matrix{u&c_1\cr c_2&v}\right)$ of $M_3(\ce)$, where $c_1$ and
$c_2$ are compact. It follows that $\left(\matrix{u&0\cr
0&v\cr}\right)$ is Fredholm, so $u$ and $v$ are Fredholm liftings of
the images in $\co$ of $U$ and $V$. But this contradicts the last part
of Lemma 21.  \hfill $\square$ \bigskip

The proof of Theorem 26 generalizes in a straightforward way to give,
for every $k\in\N$, an example of a space $X$ such that $X^n$ is
isomorphic to $X^m$ if and only if $m=n$ (mod $k$). It is likely that
every Fredholm operator on the space $X$ of this section has zero
index, so that $X$ is not isomorphic to its hyperplanes. Working with
a dyadic tree may then give an example of a space $X$ isomorphic to
$X^2$ but not isomorphic to its hyperplanes.

To end this subsection, we explain, as promised, how K-theory can be
used to prove Theorem 26. (This argument also generalizes easily to
deal with the spaces mentioned in the last paragraph.) If $A$ is a
unital Banach algebra and $e$ an idempotent in $M_n(A)$ for some $n\ge
1$, let $[e]$ denote the image of $e$ in the additive group $K_0(A)$;
in particular, let $[1_A]$ or simply $[1]$ denote the image in
$K_0(A)$ of the unit of $A$. We work with complex scalars for the rest
of this section.

Now suppose that Theorem 26 is false, let $U:X\rightarrow X\oplus X$
be an onto isomorphism and let $V: X\oplus X\rightarrow X$ be its
inverse. In $M_2(L(X))$ we have the equations 
$$ e = \left(\matrix{ 1 & 0\cr 0 &
1\cr}\right) = \left(\matrix{U & 0\cr}\right)
\left(\matrix{V\cr0\cr}\right); \ f = \left(\matrix{ 1 & 0\cr 0 &
0\cr}\right) = \left(\matrix{V\cr 0\cr}\right) \left(\matrix{U &
0\cr}\right).$$ 

This means that the two idempotents $e$ and $f$ in $M_2(L(X))$ are
equivalent, and this implies by definition of the addition in $K_0$
that $[1]+[1] = [1]$, so $[1] = 0$ in $K_0(L(X))$.  Taking the image
under $\theta\circ\phi:L(X)\ra\co$, this yields $[1_{\cal O}] = 0 $ in
$K_0({\cal O})$.  All we have to show now is that $[1_{\cal O}] \ne
0$.

For this we follow the proof given by Cuntz for Theorem 3.7 of [C].
(For the K-theory details we assume, see for example [B].)
By the definition of equivalence for idempotents, $1_{\cal E} = S_i^* 
S_i$ and $S_i S_i^*$ are equivalent. The relation $I - P = \sum_{i=0}^2 
S_i S_i^*$ implies in $K_0({\cal E})$ that
$$ [1_{\cal E}]-[P] = 3 [1_{\cal E}],$$
and therefore that $[P] = -2[1_{\cal E}]$.
Now consider the short exact sequence
$$ 0\rightarrow {\cal I}\fl{j} {\cal E} \fl{\pi} {\cal O}\rightarrow 
0$$
and the corresponding exact sequence in K-theory 
$$ K_1({\cal O}) \fl{\partial_1} K_0({\cal I}) \fl{j_*} K_0({\cal E})
\fl{\pi_*} K_0({\cal O}) \fl{\partial_0} K_1({\cal I}).$$
It is easy to see that $K_1({\cal I}) = 0$ and $K_0({\cal I})\simeq \Z$ 
as they are for the ideal of compact operators on $\ell_2$.
Corollary 24 and the definition of $\partial_1$ immediately
imply that $\partial_1 =0$, so 
we get an exact sequence
$$ 0\rightarrow K_0({\cal I}) \fl{j_*} K_0({\cal E})
\fl{\pi_*} K_0({\cal O}) \rightarrow 0.$$
Now, $r = [P]$ generates $j_*(K_0({\cal I})) = \ker\pi_* \simeq \Z$.
If $ 0 = [1_{\cal O}] = 
\pi_*([1_{\cal E}])$, it follows by exactness that $[1_{\cal E}] = 
nr$ for some integer $n\in\Z$. But we know that $r = -2[1_{\cal 
E}]$, so $(2n+1)r =0$, contradicting the fact that
$r$ generates a group isomorphic to $\Z$.

\bigskip

\noindent {\bf \S 5. Further results, remarks and questions}
\smallskip

\noindent {\bf (5.1)} We begin by considering a space defined in [G1],
which has an unconditional basis but fails to be isomorphic to any
proper subspace. Here, we prove the stronger result that every
operator on the space is the sum of a diagonal operator and a strictly
singular one. Notice that this result is very much in the same spirit
as Theorem 5. Consider the (non-proper) set $\cs$ of all diagonal
operators with $\pm 1$ entries. Let $X$ be the space $X(\cs)$. There
is no problem about this as the conditions not satisfied by this $\cs$
were not needed in the definition of $X(\cs)$, but rather in the proof
of Theorem 5. It is easy to see that $X$ has a 1-unconditional
basis. Now let $T$ be any operator on $X$ with zeros down the
diagonal. We shall show that $T$ is strictly singular, which then
shows that $U-\diag(U)$ is strictly singular for any $U$.

Before we state the next lemma, we remind the reader of our convention
that if $A$ is a subset of $\N$, then $A$ also denotes the projection
$\sm n \infty a_n \bfe_n\mapsto\sum_{n\in A}a_n \bfe_n$.

\proclaim Lemma 27. Let $\sq x\in\cl$ be a sequence of successive
vectors. Let $A_n=\supp(x_n)$ and, for each $n$, let $B_n\cup C_n$ be
a partition of $A_n$ into two subsets.  Then $C_nTB_nx_n\ra 0$.

\Proof  If this is not true, then we can assume, after passing to a
subsequence, that $\nm{C_nTB_nx_n}>\e$ for some fixed $\e>0$ and
every $n$. We may also perturb $T$ by a strictly singular amount so
that the matrix of $T$ has finitely many non-zero entries in each row
and column, $T$ still has zeros down the diagonal and the above 
inequality still holds. We may now pass to a further subsequence
such that all $x_n$ are disjointly supported and so are all $Tx_n$.

Now let $y_n=B_nx_n$ and $z_n=C_nTB_nx_n$ for every $n$. Then $\nrm
{y_n}{n}\le 1$, $\nm{z_n}\ge\e$ and $y_n$ and $z_n$ are disjointly
supported. Let $B=\bigcup B_n$ and $C=\bigcup C_n$ and let $U$ be the
operator $CTB$, so that $U(y_n)=z_n$ for every $n$. We now construct
a special sequence in what is becoming the usual way. (See for example 
the proof of Lemma 8.)

Given $N\in L$, an $N$-{\sl pair} is a pair $(w,w^*)$ constructed as
follows. Let $y_{n_1},\dots,y_{n_N}$ be a subsequence of $\sq y$
satisfying the RIS(1) condition. Let
$w=N^{-1}f(N)(y_{n_1}+\dots+y_{n_N})$ and let
$w^*=f(N)^{-1}(z_{n_1}^*+\dots+z_{n_N}^*)$, where each $z_{n_i}^*$ is
a support functional for $z_{n_i}$ with
$\supp(z_{n_i}^*)\subset\supp(z_{n_i})$. In this case, we have $w^*\in
A_N^*(X)$. By Lemma 7, we also have $\nrm w {\sqrt N}\le 8$. Notice
also that $w^*(Uw)=N^{-1}\sm i Nz_{n_i}^*(z_{n_i})>\e$
and $|w^*|(|w|) = 0$.

Now for any $k\in K$ we can choose a sequence $((w_i,w_i^*))_{i=1}^k$
of such pairs as follows. Let $N_1 =j_{2k}\in L$
and let $(w_1,w_1^*)$ be an
$N_1$-pair.
We know that $\sqrt{N_1}>(2k/f'(1))2^{k^2}$.
Perturb $w_1^*$ slightly so that it is in $\bf Q$, it has
the same support, still satisfies $w_1^*(U w_1)>\e$ and so that, setting
$N_2=\sigma(w_1^*)$, we have $\sqrt{f(N_2)}>2|\ran(w_1)|$. Once the
first $i-1$ pairs have been chosen, let
$N_i=\sigma(w_1^*,\dots,w_{i-1}^*)$, let $(w_i,w_i^*)$ be an
$N_i$-pair with $w_i>w_{i-1}$ and let $w_i^*$ be perturbed so that
$w_i^*\in\bf Q$, the support is the same, $w_i^*(U w_i)>\e$ and so that
$\sqrt{f(N_{i+1})}$ will be at least $2|\ran(w_1+\dots+w_i)|$. This
construction guarantees that $(w_1^*,\dots,w_k^*)$ is a special
sequence and that ${1\over 8}(\seq w k)$ satisfies the RIS(1)
condition.

By the definition of the norm, we then get that 
$$\Bnm{\sm i kUw_i}\ge f(k)^{-1/2}\sm i kw_i^*(Uw_i)>\e kf(k)^{-1/2}
\ .$$
On the other hand, we can let $g$ be the function given by Lemma 6 in
the case $K_0=K\setminus\{k\}$ and apply Lemma 3 to get an upper bound
for $\nm{\speq w k}$ of $24kf(k)^{-1}$. This shows that $\nm
U>(\e/24)f(k)^{1/2}$ for every $k\in K$, contradicting the continuity
of $T$.

To do this, let $w=(1/8)(\speq w k)$. By Lemma 3, it is enough to show
that, given any special functional
$u^*=f(k)^{-1/2}(u_1^*+\dots+u_k^*)$ of size $k$, and any interval
$E$, we have $|u^*|(|Ew|)<1$. Let $t$ be maximal such that
$u_t^*=w_t^*$. Then, for $i\le t$ we have, by the conditions on the
supports of the $y_n$ and $z_n$, that $u_i^*(w_j)=0$ for every $j$. If
$i>t+1$, then $u_i^*\in A_N^*$ for some $N\in L$ distinct from all of
$\seq N k$. This gives a good upper bound for $|u_i^*|(|w_j|)$. If
$N<N_j$ then $N<\sqrt{N_j}$ by choice of $L$. Since
$\nrm{\,|w_j|\,}{\sqrt{N_j}}\le 8$, we have $|u_i^*|(|w_j|)\le 8
f(N)^{-1}$. If $N>N_j$ then $N>>N_j$ and, by Lemma 1,
$|u_i^*|(|w_j|)\le 3 f(N_j)^{-1}$. These bounds are certainly good
enough to give $|u^*|(|Ew|)\le 1$ for every such $u^*$ and finish the
proof. \hfill $\square$ \bigskip

\proclaim Corollary 28. If $\sq x\in\cl$ then $Tx_n\ra 0$.

\Proof Let $A_n=\supp(x_n)$ as before. Suppose $|A_n|$ is even. Then
because $\diag(T)=0$, we have that $Tx_n$ is four times the average of
all vectors of the form $C_nTB_nx_n$, where $B_n\cup C_n=A_n$ and
$|B_n|=|C_n|=|A_n|/2$. It follows that $\nm{Tx_n}$ is at most four
times the maximum of $\nm{C_nTB_nx_n}$, which converges to zero by
Lemma 27. If $|A_n|$ is odd, only a small modification is needed. One
can average over $C_nTB_nx_n$ where $B_n\cup C_n$ is a partition of
$A_n$ and $|B_n|$ and $|C_n|$ differ by at most 1. Then four above
must be replaced by $4n^2/(n^2-1)$. \hfill $\square$ \bigskip

\proclaim Theorem 29. Let $T$ be any operator on the space $X$.
Then $T-\diag T$ is strictly singular. 

\Proof  The previous corollary shows that $(T-\diag T)x_n\ra 0$
for every $\sq x\in\cl$. But, by Lemma 4, this implies that 
$T-\diag T$ is strictly singular. \hfill $\square$ \bigskip

\proclaim Corollary 30. The space $X$ is not isomorphic to any
proper subspace of itself.

\Proof Let $T$ be an isomorphism on to some subspace of $X$. If zero
were infinitely singular for $\diag T$, then it would also be for $T$,
by Theorem 29. Hence the diagonal entries of the matrix of $T$ are
eventually bounded below in modulus. It follows that $\diag T$ is
Fredholm. Moreover, both the rank and corank of $\diag T$ are just the
number of zeros on the diagonal, so it has index zero.  Therefore, by
Theorem 29, $T$ is Fredholm with index zero, so the subspace must be
the whole of $X$. \hfill $\square$ \bigskip

\noindent {\bf (5.2)} In this subsection we suggest possibilities for
further work along the lines of this paper and we make some
remarks. Some people may object to our new prime Banach space on the
grounds that it has no non-trivial complemented subspaces. Being prime
under these circumstances is not such a great achievement. One
possible answer to this objection is to use Theorem 5 for a slightly
larger proper set than the one we used in (4.2). There are various
possibilities. One promising one is to let $\cs$ be the set of all
spreads $S_{A,B}$, where $A$ and $B$ are infinite arithmetic
progressions in $\N$. We do not have a proof that this space is prime,
but it looks likely.

The result of (5.1) is sufficiently similar to the applications of
Theorem 5 to suggest that Theorem 5 can be generalized. It would be
nice, for example, to avoid the technical restriction on proper sets
(to do with pairs of integers). It can certainly be relaxed somewhat,
but this was not necessary for us. A motivation for carrying out such
a generalization is that it ought then to be possible to give further
interesting examples of Banach spaces with an unconditional basis. For
example, suppose one takes $\cs$ to be the semigroup generated by the
proper set from (4.4) and all diagonal maps with $\pm 1$ down the
diagonal. It seems likely that the resulting space $X(\cs)$, a sort of
combination of the spaces from (4.4) and (5.1), would be a space with
an unconditional basis isomorphic to its cube but not its square.

Such a space together with its square would be the first example of a
solution to the Schroeder-Bernstein problem using spaces with an
unconditional basis.  Recall that the Schroeder-Bernstein problem for
Banach spaces asks whether, if $X$ and $Y$ are two Banach spaces
isomorphic to complemented subspaces of each other, they must be
isomorphic. A counterexample was given in [G2]. The spaces in (4.3)
and (4.4) also produce counterexamples.

Another possibility would be to take $\cs$ to be the semigroup
generated by all spreads and all $\pm 1$-diagonal matrices. This might
well be a prime Banach space with an unconditional basis, which 
would surely be a genuine example by any standards.

Going back to proper sets of spreads, note that for any such set
$\cs$, the space $X(\cs)$ has no unconditional basic sequence. To see
this, suppose that $Y$ is a subspace generated by an unconditional
block basis. Then $L(Y,X)$ is not separable, even in the $\nmm.$-norm.
Since $\cs$ is countable, this contradicts Theorem 5. It follows from
[G3] that $X(\cs)$ has a hereditarily indecomposable subspace. In
other words, the extra structure given to these spaces by $\cs$
disappears when one passes to an appropriate subspace. It can in
fact be shown that for a proper set $\cs$ the space $X(\cs)$ has
a hereditarily indecomposable subspace generated by a subsequence
of the unit vector basis.
\bigskip

\noindent {\bf References}
\def \item{\noindent}

\item {[B]}  B. Blackadar, {\sl K-theory for operator algebras},
MSRI Publications 5, Springer Verlag 1986.

\item {[C]} J. Cuntz, {\sl K-theory for certain ${\rm C}^*$-algebras},
Ann. of Math. {\bf 113} (1981), 181-197.

\item {[G1]} W. T. Gowers, {\sl A solution to Banach's hyperplane 
problem}, Bull. L.M.S., to appear.

\item {[G2]} W. T. Gowers, {\sl A solution to the Schroeder-Bernstein
problem for Banach spaces}, submitted.

\item {[G3]} W. T. Gowers, {\sl A new dichotomy for Banach spaces}, 
preprint.

\item {[GM]} W. T. Gowers and B. Maurey, {\sl The unconditional basic
sequence problem}, Journal A.M.S. {\bf 6} (1993), 851-874.

\item {[L]} J. Lindenstrauss, {\sl On complemented subspaces of $m$},
Israel J. Math. {\bf 9} (1971), 279-284.

\item {[LT]} J. Lindenstrauss and L. Tzafriri, {\sl Classical Banach 
spaces I: sequence spaces}, Springer Verlag (1977).

\item {[P]} A. Pe\l czy\'nski, {\sl Projections in certain Banach spaces},
Studia Math. 19 (1960), 209-228.

\item {[S]} T. Schlumprecht, {\sl An arbitrarily distortable Banach 
space}, Israel J. Math. {\bf 39} (1991), 81-95.

\bye